\documentclass[useAMS]{gCOM2e}

\usepackage{algorithm}

\begin{document}
\doi{0020716YYxxxxxxxx}
 \issn{1029-0265}
\issnp{0020-7160} \jvol{00} \jnum{00} \jyear{2011} \jmonth{July}

\markboth{Ramakrishna Tipireddy, Eric T. Phipps and Roger G. Ghanem}{International Journal of Computer Mathematics}


\title{{\itshape Solvers and precondtioners based on Gauss-Seidel and Jacobi algorithms for non-symmetric stochastic Galerkin system of equations}}

\author{Ramakrishna Tipireddy$^{\rm a,\ast}$ \thanks{$^\ast$Corresponding author. Address: 3620 S. Vermont Avenue, KAP 210, Los Angeles, CA 90089-2531. Email: tipiredd@usc.edu Tel: (213) 740-9615
\vspace{6pt}}, Eric T. Phipps$^{\rm b}$ and Roger G. Ghanem$^{\rm a}$ \\\vspace{6pt}  $^{\rm a}${\em{Department of civil and environmental engineering, University of Southern California, Los Angeles, CA, USA}}; \\\vspace{3pt} $^{\rm b}${\em{Sandia National Laboratories$^\dag$, Albuquerque, NM, USA\thanks{$^\dag$Sandia National Laboratories is a multi-program laboratory managed and operated by Sandia Corporation, a wholly owned subsidiary of Lockheed Martin Corporation, for the U.S. Department of Energy's National Nuclear Security Administration under contract DE-AC04-94AL85000.}}}\\\vspace{6pt}}

\maketitle

\clearpage
\begin{abstract}


In this work, solvers and preconditioners based on Gauss-Seidel and Jacobi algorithms are explored for stochastic Galerkin discretization of  partial differential equations (PDEs) with random input data. Gauss-Seidel and Jacobi algorithms are formulated such that the existing software is leveraged in the computational effort. These algorithms are also used as preconditioners to Krylov iterative methods. The solvers and preconditiners are compared with Krylov based iterative methods with the traditional mean-based preconditioner~\cite{RK:Powell2009} and Kronecker-product preconditioner~\cite{RK:Ullmann2010} by solving a steady state state advection-diffusion equation, which upon discretization, results in a non-symmetric positive definite matrix on left-hand-side. Numerical results show that an approximate version of Gauss-Seidel algorithm is a good preconditioner for GMRES to solve non-symmetric Galerkin system of equations.


\end{abstract}

\section{Introduction} \label{RK:sec:intro}

To account for uncertainties in model, data and parameters, physical phenomena are often modeled as partial differential equations (PDEs) with random coefficients.  While Monte Carlo sampling methods remain the standard for providing a probabilistic characterization to the solution of these problems, there is growing need for approaches that address the computational challenge associated with statistical sampling of large scale computational models.  More recently, methods based on stochastic Galerkin approximations~\cite{RK:Ghanem1991, RK:Babuska2004} have been increasingly explored as an alternative to Monte Carlo for solving these problems because of their computational and analytical properties, as well as for their ability to provide an approximation to the deterministic mapping that exists between random input parameters and random solution. 

Stochastic finite element methods based on intrusive stochastic Galerkin methods~(\cite{RK:Rosseel2010, RK:Pellissetti2000, RK:Eiermann2007, RK:Ghanem1996}) and non-intrusive stochastic collocation methods~(\cite{RK:Smolyak1963, RK:Xiu2005, RK:Nobile2008, RK:Babuska2007}) have gained popularity in recent years.  Here the intrusive and the non-intrusive methods are defined as methods which require significant changes to software already available for solving deterministic PDEs and methods that can use this legacy software without modifications, respectively. Both methods exploit solution regularity to improve on the convergence rates of Monte Carlo methods.  The first approach translates the stochastic PDE into a coupled set of deterministic PDEs while the second samples the stochastic PDE at a predetermined set of collocation points, resulting in a set of uncoupled deterministic PDEs. The solution at these collocation points is then used to interpolate the solution in the entire random input domain. Extending legacy software to support stochastic collocation methods is simpler than supporting stochastic Galerkin methods (SGMs). Moreover, intrusive SGMs require specialized linear solvers. However, the resulting set of PDEs in the stochastic Galerkin system is much smaller in number than that in the collocation method. For a canonical random diffusion problem, it is shown~\cite{RK:Elman2011} that SGM using iterative Krylov-based linear solvers and mean-based preconditioning~\cite{RK:Powell2009} is more efficient than the non-intrusive sparse grid collocation method.    

While the stochastic Galerkin method is often considered to be a fully intrusive method, there are in fact a variety of solver approaches for the stochastic Galerkin method with varying level of intrusiveness.  In the present work, the less intrusive Gauss-Seidel and Jacobi mean-based algorithms are explored to solve stochastic Galerkin system of equations.  We consider these methods to be less intrusive than the Krylov-based methods as they allow reuse of existing deterministic solvers.  Moreover preconditioning techniques for Krylov-based methods based on Gauss-Seidel and Jacobi ideas are also explored and compared to traditional mean-based preconditioning.  Various iterative methods and preconditioners based on matrix splitting methods to solve stochastic Galerkin system of equations are compared in~\cite{RK:Rosseel2010}. The $C_i$-splitting methods discussed in section~3.1.3 in~\cite{RK:Rosseel2010} are similar to the Gauss-Seidel and Jacobi algorithms discussed in this paper. In~\cite{RK:Rosseel2010}, it was reported that the mean-based preconditioner was more efficient in terms of computational time than the other $C_i$-splitting methods considered.  A symmetric Gauss-Seidel preconditioner was used in~\cite{RK:Rosseel2010} which involves one forward and one backward Gauss-Seidel iteration. In the present paper, we use only one forward Gauss-Seidel iteration for preconditioning and hence the computational cost of this preconditioner is half that of the symmetric Gauss-Seidel preconditioner.   Although, these solvers and preconditioners can be used to solve any stochastic Galerkin system of equations, we focus here on solving non-symmetric stochastic Galerkin system of equations. Gauss-Seidel and Jacobi algorithms are carefully formulated such that the number of matrix-vector products are minimized. Unlike the mean-based preconditioner, which has information from the mean stiffness matrix only, preconditioners based on Gauss-Seidel and Jacobi algorithms make use of higher order information of the stochastic stiffness matrix. Another preconditioner that makes use of higher order information is the Kronecker product preconditioner defined in~\cite{RK:Ullmann2010}.  This preconditioner has the Kronecker product structure of the Galerkin stiffness matrix and is constructed to be close to the stochastic stiffness matrix.  All these techniques are then compared by solving a non-symmetric test problem with random diffusion coefficient.  These comparisons demonstrate a trade-off in computational cost versus intrusiveness with the Krylov-based methods using an approximate Gauss-Seidel or Jacobi mean preconditioner being the most efficient. Main contribution in this paper is to adapt Gauss-Seidel and Jacobi algorithms such that the legacy software can be reused to solve stochastic partial differential equations and to reduce the computational cost of Krylov based iterative methods by using these algorithms as preconditioners in which the number of matrix-vector products is reduced by minimizing the duplicated computations. 

This paper is organized as follows.  Two models of the input random field are developed in section~\ref{RK:eq:rf} which dictate very different behavior for the stochastic solution methods considered next.  Section~\ref{RK:sec:sgm} describes the stochastic Galerkin method. Various solver and preconditioning methods for solving the stochastic Galerkin system of equations are introduced in section~\ref{RK:sec:sgmsoln}.  In section~\ref{RK:sec:probstate}, a test problem governed by the advection-diffusion equation with a random diffusion coefficient is formulated.  In section~\ref{RK:sec:numer}, numerical experiments are carried out to compare the efficiency of various solver and preconditioning methods that have been introduced in section~\ref{RK:sec:sgmsoln}. Finally section~\ref{RK:sec:conc} provides the concluding remarks. 
 
\section{Stochastic partial differential equations} \label{RK:sec:spde}

Let  $D$ be an open subset of $\mathbb{R}^n$ 
and $(\Omega, \Sigma, P)$ be a complete probability space with sample space $\Omega$, $\sigma$-algebra $\Sigma$ and  probability measure $P$. We are interested in studying the following stochastic partial differential equation: find a random field, $u(x,\omega):D \times \Omega \rightarrow \mathbb{R}$ such that the following holds $P$-almost surely ($P$-a.s.):
\begin{equation}\label{RK:eq:spdeop}
 \mathcal{L}(x,\omega;u) = f(x,\omega)  \;\; \rm{in}~D\times \Omega,
\end{equation}
subject to the boundary condition
\begin{equation}\label{RK:eq:spdebc}
 \mathcal{B}(x,\omega;u) = g(x,\omega)  \;\; \rm{on}~\partial D\times \Omega,
\end{equation}
where $\mathcal{L}$ is a differential operator and $\mathcal{B}$ is a boundary operator. 

\section{Input random field model}\label{RK:eq:rf}

Uncertainty in stochastic PDEs often arises by treating coefficients in the differential operator $ \mathcal{L}(x,\omega;u)$ as random fields.  For computational purposes, each stochastic coefficient $a(x,\omega)$ must be discretized in both spatial and stochastic domains.  To this  end, it is often approximated with a truncated series expansion that separates the spatial variable $x$ from the stochastic variable $\omega$ resulting in a representation by a finite number of random variables.   In the present problem, two cases of random field models are considered. In the first case, the random field is assumed to be charcaterized by a truncated Karhunen-Lo\`eve expansion obtained by specifying the covariance function and a product uniform measure on the associated random variables.  In the second case, the random field is assumed to have a log-normal distribution, that is $a(x,\omega) = \exp{(g(x,\omega))}$ where $g(x,\omega)$ is a Gaussian random field, and is approximated by a truncated polynomial chaos expansion.

\subsection{Karhunen-Lo\`{e}ve representation}\label{RK:sec:kl}

Let $C(x_1,x_2) = E[a(x_1,\omega)a(x_2,\omega)]$ be the covariance function of the random field $a(x,\omega)$, where $E[\cdot]$ denotes  mathematical expectation. Then $a(x,\omega)$ can be approximated through its truncated Karhunen-Lo\`{e}ve (K-L)  expansion~\cite{RK:Ghanem1991} given by 
\begin{equation}\label{RK:eq:kl}
 a(x,\omega) \approx \tilde{a}(x,\xi(\omega)) = a_0(x) + \sum_{i=1}^M \sqrt{\lambda_i} a_i(x) \xi_i(\omega),
\end{equation}
where $a_0(x)$ is the mean of the random field $a(x,\omega)$  and $\{(\lambda_i,a_i(x))\}_{i\geq1}$ are solutions of the integral eigenvalue problem
\begin{equation}\label{RK:eq:inteig}
 \int_D C(x_1,x_2) a_i(x_2)dx_2 = \lambda_i a_i(x_1).
\end{equation}
The eigenvalues $\lambda_i$ are positive and non-increasing, and the eigenfunctions $a_i(x)$ are orthonormal, that is, 
\begin{equation}\label{RK:eq:ortheigf}
 \int_D a_i(x) a_j(x) = \delta_{ij},
\end{equation}
where $\delta_{ij}$ is the Kronecker delta.  
In Eq.~\ref{RK:eq:kl}, $\{\xi_i\}_{i=1}^M$ are uncorrelated random variables with zero mean. As a first test-case, the diffusion coefficient $a(x,\omega)$ is modeled as a random field with an exponential covariance function

\begin{equation}\label{RK:eq:exp_cov}
C(x_1,x_2) = \sigma^2\exp(-\|x_1-x_2\|_1/L)\ .
\end{equation}
In addition, the random variables $\xi_i(\omega)$ in the K-L expansion are assumed to be independent identically distributed with a uniform
distribution.

\subsection{Polynomial chaos representation}\label{RK:sec:pce}

In general, the K-L expansion results in a representation that, while linear in the KL random variables, involves variables that are statistically dependent with a joint distribution that is data-dependent.  For problems driven by experimental evidence, this could present both  theoretical and algorithmic challenges.  An alternative representation relies on representing the random field directly in the form of a polynomial chaos decomposition~\cite{RK:Ghanem1991}.  Accordingly, the random field $a(x,\omega)$ is approximated as

\begin{equation}\label{RK:eq:pce}
 a(x,\omega) \approx \tilde{a}(x,\boldsymbol{\xi}(\omega)) = a_0(x) + \sum_{i=1}^{N_{\xi}} a_i(x) \psi_i(\boldsymbol{\xi}) ,
\end{equation}
where $\{\psi_i(\boldsymbol{\xi})\}$ are orthogonal polynomials with respect to measure of the random variables $\bm\xi\in{\mathbb{R}}^M$ which are chosen, typically, to be statistically independent.  Specifically,

\begin{equation}\label{RK:eq:innprod}
\langle \psi_i , \psi_j \rangle \equiv \int_{\mathbb{R}^M} \psi_i(\boldsymbol{\xi}) \psi_j(\boldsymbol{\xi}) f_{\bm\xi} d\bm\xi = \delta_{ij}.
\end{equation}
Accordingly, and as a second test-case, the diffusion coefficient $a(x,\omega)$, is modeled as a log-normal random field~\cite{RK:Ghanem1999} where $a(x,\omega) = \exp[g(x,\omega)]$ and $g(x,\omega)$ is a Gaussian random field with exponential covariance~\eqref{RK:eq:exp_cov}. Here, $g(x,\omega)$ is approximated with a truncated K-L expansion of the form

\begin{equation}
g(x,\omega) \approx \bar{g}(x,\bm\xi) = g_0(x) + \sum_{i=1}^M \sqrt{\lambda}_ig_i(x)\xi_i(\omega)\ .
\end{equation}
In this case the random variables $\xi_i$ are standard normal random variables and thus are independent and $a(x,\omega)$ can be approximated with a truncated polynomial chaos expansion~\eqref{RK:eq:pce} where $\psi_i$ are multidimensional Hermite polynomials in gaussian random variables $\bm\xi$.  For a given total polynomial order $p$, the total number of polynomials $\{\psi_i(\boldsymbol{\xi})\}$ is $N_{\xi}+1 = \frac{(M+p)!}{M!p!}$.

\section{Stochastic Galerkin method}\label{RK:sec:sgm}

Let $H^1_0(D)$ be the subspace of the Sobolev space $H^1(D)$ that vanishes on the boundary $\partial D$ and is equipped with the norm  $\|u\|_{H^1_0(D)}=[\int_D |\nabla u|^2 dx]^{\frac{1}{2}}$.   Problem~\eqref{RK:eq:spdeop} can then be written in the following equivalent variational form~\cite{RK:Ghanem2006}:  find $u \in H^1_0(D) \otimes L_2(\Omega)$ such that 
\begin{equation}\label{RK:eq:svfop}
b(u,v) = l(v), \quad \forall v \in H^1_0(D) \otimes L_2(\Omega),
\end{equation}
where $b(u,v)$ is a continuous and coercive bilinear form and $l(v)$ is a continuous bounded linear functional.
In the stochastic Galerkin method, we seek the solution of the variational problem~\eqref{RK:eq:svfop} in a tensor product space $X_h \otimes Y_p$, where, $X_h\subset H^1_0(D)$ is finite dimensional space of continuous polynomials corresponding to the spatial discretization of $D$ and $Y_p\subset L_2(\Omega)$ is the space of random variables spanned by polynomial chaos~\cite{RK:Ghanem1991} of order up to $p$. Then the finite dimensional approximation $u_{X_hY_p}(x,\omega)$ of the exact solution $u(x,\omega)$ on the tensor product space $X_h \otimes Y_p$ is given as the solution to
\begin{equation}\label{RK:eq:sgm}
b(u_{X_hY_p},v) = l(v) \quad \forall v\in X_h \otimes Y_p.
\end{equation}

In equation~\eqref{RK:eq:sgm} the input random field $a(x,\omega)$ in the bilinear form $b(u_{X_hY_p},v)$ can be approximated using either a K-L expansion or a polynomial chaos expansion depending on a choice of model for the random field.  The finite dimensional approximation $u_{X_hY_p}(x,\omega)$ is represented with truncated polynomial chaos expansion, where the multidimensional polynomial chaos are orthogonal with respect to the probability measure of the underlying random variables. The resulting set of coupled PDEs are then discretized using standard techniques such as the finite element or finite difference methods. In the present problem, the set of coupled PDEs are discretized using the finite element method and the resulting system of linear equations can be written as

\begin{equation}\label{RK:eq:disceqns}
\sum_{j=0}^{N_{\xi}}\sum_{i=0}^{\hat{P}} c_{ijk} K_i u_j =  f_k,  \quad k = 0,\cdots, N_{\xi},
\end{equation}
where $f_k = E\{f(x,\boldsymbol{\xi}) \psi_k\}$, $c_{ijk}= E\{\xi_i \psi_j \psi_k \}$ and $\hat{P} = M$ when $a(x,\omega)$ is approximated by a truncated K-L expansion, or $c_{ijk}= E\{\psi_i \psi_j \psi_k \}$ and $\hat{P} = \hat{N}_{\xi}$ when $a(x,\omega)$ is approximated by a polynomial chaos expansion. Here $\{K_i\in\mathbb{R}^{N_x\times N_x}\}_{i=0}^{\hat{P}}$ are the polynomial chaos coefficients of the stiffness matrix~(section (3.4) of \cite{RK:Powell2009}) and $\{u_j\in\mathbb{R}^{N_x}\}_{j=0}^{N_\xi}$ are the polynomial chaos coefficients of the discrete solution vector
\begin{equation}
  u_j = [ u_{0j}, \dots, u_{N_xj}]^T, \quad j=0,\dots,N_\xi.
\end{equation}
 
Equation~\eqref{RK:eq:disceqns} can be written in the form of a global stochastic stiffness matrix of size $((N_{\xi}+1) \times N_x)$ by $((N_{\xi}+1) \times N_x)$ as
\begin{equation}\label{RK:eq:stochstiff}
  \begin{bmatrix}
    K^{0,0} & K^{0,1} & \cdots & K^{0,N_{\xi}} \\
    K^{1,0} & K^{1,1} & \cdots & K^{1,N_{\xi}} \\
    \vdots & \vdots & \vdots & \vdots \\
    K^{N_{\xi},0} & K^{N_{\xi},1} & \cdots & K^{N_{\xi},N_{\xi}} \\
  \end{bmatrix}
  \times \left\lbrace\begin{array}{c} u_1 \\ u_2 \\ \vdots \\ u_{N_{\xi}} \end{array}\right\rbrace
= \left\lbrace\begin{array}{c} f_1 \\ f_2 \\ \vdots \\ f_{N_{\xi}}\end{array}\right\rbrace
\end{equation}
where $K^{j,k} = \sum_{i=0}^{\hat{P}} c_{ijk} K_i$.   We will denote this system as $\bar{K}\bar{u} = \bar{f}$.  In practice it is prohibitive to assemble and store the global stochastic stiffness matrix in this form, rather each block of the stochastic stiffness matrix can be computed from the $\{K_i\}$ when needed. The stochastic stiffness matrix is block sparse, which means that some of the off-diagonal blocks are zero matrices, because of the fact that the $c_{ijk}$ defined above vanishes for certain combination of $i, j$ and $k$. 
Unlike the off-diagonal blocks, the diagonal blocks have contributions from the mean striffness matrix $K_{0}$ and thus dominate over the off-diagonals. These properties of block sparsity and diagonal dominance are exploited to develop solvers and preconditioners based on Jacobi and Gauss-Seidel algorithms.  

\section{Solution methods for stochastic Galerkin systems}\label{RK:sec:sgmsoln}

In this section, various solver techniques and preconditioning methods for solving the linear algebraic equations arising from stochastic Galerkin discretizations~\eqref{RK:eq:disceqns} are described. The solver methods discussed are: relaxation methods, namely, a Jacobi mean method and a Gauss-Seidel mean method,  and Krylov-based iterative methods~\cite{RK:Saad1996}.  Also various stochastic preconditioners used to accelerate convergence of the Krylov methods are discussed, including mean-based~\cite{RK:Powell2009}, Gauss-Seidel mean, approximate Gauss-Seidel mean, approximate Jacobi mean and Kronecker product~\cite{RK:Ullmann2010} preconditioners. The relaxation schemes can be viewed as fixed point iterations on a preconditioned system~\cite{RK:Saad1996}. In Jacobi and Gauss-Seidel methods, mean splitting is used rather than traditional diagonal block splitting as it allows use of the same mean matrix $K_0$ for all inner deterministic solves (and thus reuse of the preconditioner $P_0 \approx K_0$). 

\subsection{Jacobi mean algorithm}\label{RK:sec:jacobi}

In this method, systems of equations of size equal to that of the deterministic system are solved iteratively by updating the right-hand-side to obtain the solution to the stochastic Galerkin system of equations~\eqref{RK:eq:disceqns}:
\begin{equation}\label{RK:eq:sfem_jacobi}
c_{kk0 }K_0 u^{new}_k = f_k - \sum_{j=0}^{N_{\xi}} \sum_{i=1}^{\hat{P}} c_{ijk} K_i u^{old}_j, \quad k = 0,\cdots,N_{\xi}.
\end{equation}
The above system of equations is solved for $k = 0,\cdots,N_{\xi}$ using any solution technique appropriate for the mean matrix $K_0$.  Thus existing legacy software can be used with minimal modification to solve the stochastic Galerkin system.  In this work, Krylov-based iterative methods with appropriate preconditioners will be used. One cycle of solves from $k = 0,\cdots,N_{\xi}$ is considered one Jacobi outer iteration, and after each outer iteration, the right-hand-side in equation~\eqref{RK:eq:sfem_jacobi} is updated replacing $\{u^{old}_j\}$ with the new solution $\{u^{new}_j\}$. These outer iterations are continued until the required convergence tolerance is achieved.  
The Jacobi mean algorithm is shown in Algorithm~\ref{RK:alg:jacobi}.
\begin{algorithm}[H]
\caption{Jacobi mean algorithm} \label{RK:alg:jacobi} 
  1. Choose initial guess $\bar{u}^{0}$ and compute residual $\bar{r}=\bar{K}\bar{u}^{0}-\bar{f}$ \\
  2. Iteration count, $itr = 0$\\
  3. {\bf while} $ \frac{\|\bar{r}\|_2}{\|\bar{f}\|_2} > {\tt tol}$ {\bf do} \\
  4.  $\quad$ {\bf for} $ k=0 \ldots N_{\xi}$ {\bf do} \\  
  5. $\quad$ $\quad$ Solve $c_{kk0} K_0 u^{(itr+1)}_k = f_k - \sum_{j=0}^{N_{\xi}} \sum_{i=1}^{\hat{P}} c_{ijk} K_i u^{(itr)}_j$ \\  
  6. $\quad$ {\bf end for}\\
  7. $\quad itr = itr + 1$\\   
  8. $\quad \bar{r}=\bar{K}\bar{u}^{itr}-\bar{f}$ \\
  9. {\bf end while}.
\end{algorithm}
Note that for a given outer iteration, all of the right-hand-sides for $k = 0,\cdots,N_{\xi}$ are available simultaneously, and thus their solution can be efficiently parallelized.  Moreover block algorithms optimized for multiple right-hand-sides may be used to further increase performance.   Finally this approach does not require a large amount of memory to compute the solution.  The disadvantage of the method is that it may not converge or may converge very slowly when the diagonal blocks of the stochastic stiffness matrix are less dominant over off-diagonal blocks. 

\subsection{Gauss-Seidel mean iterative method}\label{RK:sec:gs}

The Gauss-Seidel method considered is similar to the Jacobi method, except the right-hand-side in equation~\eqref{RK:eq:sfem_jacobi} is updated after each deterministic solve with the newly computed $u^{new}_k$. Symbolically this is written
\begin{equation}\label{RK:eq:sfem_gs}
c_{kk0 }K_0 u^{new}_k = f_k - \sum_{j=0}^{k-1} \sum_{i=1}^{\hat{P}} c_{ijk} K_i u^{new}_j - \sum_{j=k}^{N_{\xi}} \sum_{i=1}^{\hat{P}} c_{ijk} K_i u^{old}_j, \quad k = 0,\cdots,N_{\xi}.
\end{equation}
As before, one cycle of solves from $k = 0,\cdots,N_{\xi}$ is considered one outer iteration of the Gauss-Seidel method, and these outer iterations are repeated until the required convergence tolerance is achieved. Note however that computing the updates as shown here would result in a large number of duplicated matrix-vector products $K_i u_j$ for each outer iteration.  Instead, after each $u_k^{new}$ is computed by solving the mean linear system, we first compute $y = K_i u_k^{new}$ for all $i$ in which $c_{ijk}$ is nonzero for any $j$.  Then for each corresponding $j$ we update $f_j \leftarrow f_j - c_{ijk} y$.  This allows all of the right-hand-sides to be updated as required using the fewest number of matrix-vector products and without resorting to storing intermediate products.
The complete Gauss-Seidel algorithm is shown in Algorithm~\ref{RK:alg:gs}.
\begin{algorithm}[H]   
\caption{Gauss-Seidel mean algorithm} \label{RK:alg:gs} 
  1. Choose initial guess $\bar{u}^{0}$ and compute residual $\bar{r} = \bar{f} - \bar{K} \bar{u}^{0}$ \\
  2. Iteration count, $itr = 0$\\
  3. Initialize $\bar{z} = \bar{r}$ \\
  4. {\bf while} $ \frac{\|\bar{r}\|_2}{\|\bar{f}\|_2} > {\tt tol}$ {\bf do} \\
  5. $\quad$ $\bar{r} = \bar{f}$ \\
  6. $\quad$ {\bf for} $ k=0 \ldots N_{\xi}$ {\bf do} \\  
  7. $\quad\quad$ Solve $c_{kk0} K_0 u_k = z_k$ \\  
  8. $\quad\quad$ $z_k = f_k$ \\
  9. $\quad\quad$ {\bf for} $ i = 1,\dots,\hat{P}$ \\
  10. $\quad\quad\quad$ $y = K_i u_k$ \\
  11. $\quad\quad\quad$ {\bf for} $ j = 1,\dots,N_\xi$ \\
  12. $\quad\quad\quad\quad$ {\bf if} $c_{ijk}\neq 0$ {\bf then} \\
  13. $\quad\quad\quad\quad\quad$ $z_j = z_j - c_{ijk} y$ \\
  14. $\quad\quad\quad\quad\quad$ $r_j = r_j - c_{ijk} y$ \\
  15. $\quad\quad\quad\quad$ {\bf endif} \\
  16. $\quad\quad\quad$ {\bf end for}\\
  17. $\quad\quad$ {\bf end for}\\
  18. $\quad\quad$ $r_k = r_k - c_{kk0} K_0 u_k$ \\
  19. $\quad$ {\bf end for}\\
  20. $\quad$ $itr = itr + 1$\\   
  21. {\bf end while}.
\end{algorithm}
Often this method converges in fewer iterations than the Jacobi method, at the expense of no longer having all of the right-hand-sides available simultaneously.  Unlike diagonal block splitting methods defined in~\cite{RK:Rosseel2010}, mean block splitting is used in both Jacobi and Gauss-Seidel algorithms and hence the left-hand-side matrix is the mean matrix for all inner deterministic problems and only the right-hand-side changes. In such cases recycled Krylov basis methods could be explored to increase performance.

\subsection{Krylov based iterative methods with matrix-free operations}\label{RK:sec:gmres}
Krylov based iterative methods~\cite{RK:Saad1996} such as the conjugate gradient (CG) method and the generalized minimal residual (GMRES) method can be used to solve the stochastic Galerkin system~\eqref{RK:eq:disceqns} in which matrix vector products $\bar{v} = \bar{K}\bar{u}$ are computed using ``matrix-free'' operations:
\begin{equation}\label{RK:eq:matfree}
v_k = \sum_{j=0}^{N_{\xi}} \sum_{i=0}^{\hat{P}} c_{ijk} K_i u_j, \quad k = 0,\cdots,N_{\xi}.
\end{equation}
If the matrix vector products are computed from Eq.~\ref{RK:eq:matfree}, it is not required to assemble the full stochastic Galerkin stiffness matrix, drastically decreasing memory requirements.  However if a large number of iterations of a Krylov method such as GMRES are required, allocation of the Krylov basis may still require a very large amount of memory.  Thus good preconditioning strategies for the stochastic Galerkin system are required, several of which will be discussed below.

\subsubsection*{Mean-based preconditioner}\label{RK:sec:mb}
The mean-based preconditioner~\cite{RK:Powell2009} is given by $P = \rm{diag} \{P_0, \cdots, P_0 \}$ where $P_0\approx K_0$ is a preconditioner for the mean.  The mean-based preconditioner is very efficient to compute and apply, since it only must be generated once from a matrix that is of the size of the deterministic system. However it doesn't incorporate any higher-order stochastic information, thus its performance degrades as the stochastic dimension, polynomial order, or random field variance increases~\cite{RK:Pellissetti2000,RK:Ullmann2010}.

\subsubsection*{Gauss-Seidel preconditioner}\label{RK:sec:prec_gs}
One or more outer iterations of the Gauss-Seidel mean algorithm can be used as a preconditioner to the Krylov based iterative methods. An advantage of this method is that the cost of applying the preconditioner can be controlled by adjusting the tolerance of the inner deterministic solves and number of outer iterations. Decreasing this tolerance and increasing the number of outer iterations will reduce the number of iterations in the Krylov method, but make the preconditioner more expensive to apply, and thus these must be balanced to minimize overall computational cost.  Generally we have found the cost of the preconditioner to be dominated by solving the mean systems, and thus the performance was improved by loosening the outer solver tolerance or limiting the number of outer iterations.  In the results presented below we limited the preconditioner to only one Gauss-Seidel iteration. 

\subsubsection*{Approximate Gauss-Seidel preconditioner}\label{RK:sec:ags}
The process of increasing the inner solver tolerance can be taken to its extreme of replacing the inner mean solves by application of the mean preconditioner.  As with the Gauss-Seidel preconditioner above, we found experimentally that this approach worked best with only one Gauss-Seidel iteration, and adding additional iterations did not improve the quality of the preconditioner.  We also found that the cost of the preconditioner was reduced dramatically if only the first-order terms in the expansion for the stiffness matrix were used in the preconditioner and using higher-order terms did not improve performance.  We refer to this as the approximate Gauss-Seidel preconditioner.  

\subsubsection*{Approximate Jacobi preconditioner}\label{RK:sec:aj}
Similar to the approximate Gauss-Seidel preconditioner, Jacobi iterations can be used with a preconditioner in place of the mean stiffness matrix.  In this case we used two outer Jacobi iterations, since the first iteration is equivalent to mean-based preconditioning (i.e., the additional terms on the right-hand-side of equation~\eqref{RK:eq:sfem_jacobi} are zero).  Increasing the number of outer iterations did not improve the efficiency of the overall solver.  We refer to this as the approximate Jacobi preconditioner. Both  approximate Gauss-Seidel and approximate Jacobi preconditioners are found to be very effective in reducing the number of Krylov iterations and are also not very expensive to apply. 

\subsubsection*{Kronecker product preconditioner}\label{RK:sec:kp}
The Kronecker product preconditioner~\cite{RK:Ullmann2010} is defined as $P_1 = G \otimes K_0$,  where $K_0$ is the mean stiffness matrix and $G$ is  
\begin{equation}
 G = \sum_{i=0}^{\hat{P}}  \frac{\text{tr}(K_i^T K_0)} {\text{tr}(K_0^T K_0)} G_i
\end{equation}
where, $G_i(j,k) = c_{ijk}$. Unlike the mean-based preconditioner, the Kronecker product preconditioner incorporates higher order stochastic information allowing the Krylov based algorithm to converge in fewer iterations. However the disadvantage is that the cost of constructing the Kronecker product preconditioner is greater than that for the mean-based preconditioner, and is also more expensive to apply. Over all, solution time is found to be less than that with mean-based preconditioner as reported in~\cite{RK:Ullmann2010}. However, we found in our numerical experiments that the approximate Gauss-Seidel and approximate Jacobi preconditioner performed better than the Kronecker product preconditioner. 

\section{Problem Statement} \label{RK:sec:probstate}
In this work, a stochastic steady state advection-diffusion equation~\cite{RK:Xu2002} with Dirichlet boundary conditions is used as a test problem for various solution methods.  Assume $a(x,\omega):D\times \Omega \rightarrow \mathbb{R}$ to be a random field that is bounded and strictly positive, that is,
\begin{equation}\label{RK:eq:rf_bound}
0 < a_l \leq a(x,\omega) \leq a_u < \infty \quad \rm{a.e.} \quad \rm{in} \quad D\times \Omega.
\end{equation}
For the steady-state advection-diffusion equation, we wish to compute a random field $u(x,\omega):D\times \Omega \rightarrow \mathbb{R}$,  such that the following holds $P$-almost surely ($P$-a.s.):
\begin{align}\label{RK:eq:spde_ad}
\vec{w}\cdot \nabla u(x,\omega) -\nabla . (a(x,\omega) \nabla u(x,\omega))&=f(x,\omega) \;\; \rm{in}~D\times \Omega, \\
u(x,\omega)&=u_0 \;\; \rm{on}~\partial D\times \Omega.
\end{align} 
where $\vec{w} = [w_x, w_y]^T$ is the advection direction, $w_x, w_y > 0$.

Problem~\eqref{RK:eq:spde_ad} can then be written in the following equivalent variational form~\cite{RK:Ghanem2006}:  find $u \in H^1_0(D) \otimes L_2(\Omega)$ such that 
\begin{equation}\label{RK:eq:svf}
b(u,v) = l(v), \quad \forall v \in H^1_0(D) \otimes L_2(\Omega),
\end{equation}
where $b(u,v)$ is the continuous and coercive (from assumption~\eqref{RK:eq:rf_bound}) bilinear form given by
\begin{equation}\label{RK:eq:bilinear_ad}
b(u,v) = E\left[\int_D\vec{w}\cdot \nabla u v dx \right] + E\left[\int_D a \nabla u \cdot \nabla v dx \right], \quad \forall u, v \in H^1_0(D) \otimes L_2(\Omega),
\end{equation}
and $l(v)$ is the continuous bounded linear functional given by 
\begin{equation}\label{RK:eq:linear}
l(v) = E\left[\int_D f v dx \right], \quad \forall v \in H^1_0(D) \otimes L_2(\Omega).
\end{equation}
From the Lax-Milgram lemma, equation~\eqref{RK:eq:svf} has unique a solution in $H^1_0(D) \otimes L_2(\Omega)$.

\section{Numerical illustration}\label{RK:sec:numer}
To compare the performance of different solvers and preconditioners discussed above, a 2-D stochastic diffusion equation and a 2-D stochastic advection-diffusion equation presented in section~\ref{RK:sec:probstate} are solved using the stochastic Galerkin method described in section~\ref{RK:sec:sgm}.  In the above two problems, the diffusion coefficient is modeled as both a random field discretized using a truncated K-L expansion with uniform random variables (section~\ref{RK:sec:kl}) and a log-normal random field discretized using a truncated polynomial chaos expansion (section~\ref{RK:sec:pce}).  In the first case, the orthogonal polynomials used in the stochastic Galerkin method are tensor products of 1-D Legendre polynomials and in the second case tensor products of Hermite polynomials are used.   For simplicity a constant unit force $f(x,\omega) = 1$ is used as the right-hand-side in equation~\eqref{RK:eq:spde_ad}.  The spatial dimensions are discretized using standard finite element mesh with linear quadrilateral elements. In advection-diffusion equation, the parameter, $\vec{w} = [1,1]^T$. The corresponding stochastic Galerkin linear system is constructed using the Stokhos and Epetra packages in Trilinos.  For the Jacobi solver, Gauss-Seidel solver and Gauss-Seidel preconditioner the non-symmetric linear systems obtained from the discretization of stochastic advection-diffusion equation are solved via multi-grid preconditioned GMRES provided by the AztecOO and ML packages in Trilinos.  

\begin{figure}[h]
\begin{center}
\subfigure[Stochastic dimension = 3]
{
\includegraphics[scale=.32]{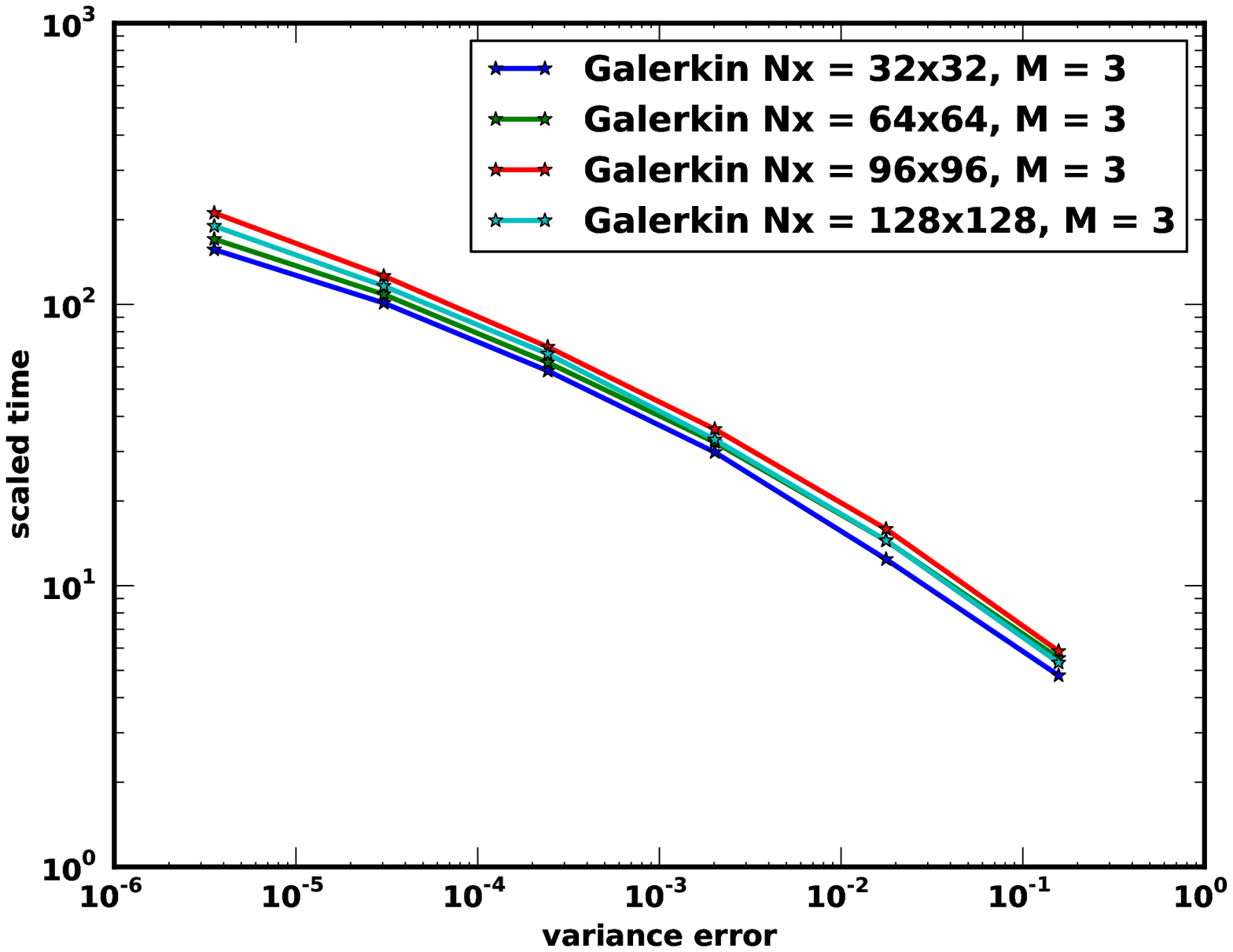}\label{RK:fig:Figure1a}
}
\subfigure[Stochastic dimension = 4]
{
\includegraphics[scale=.32]{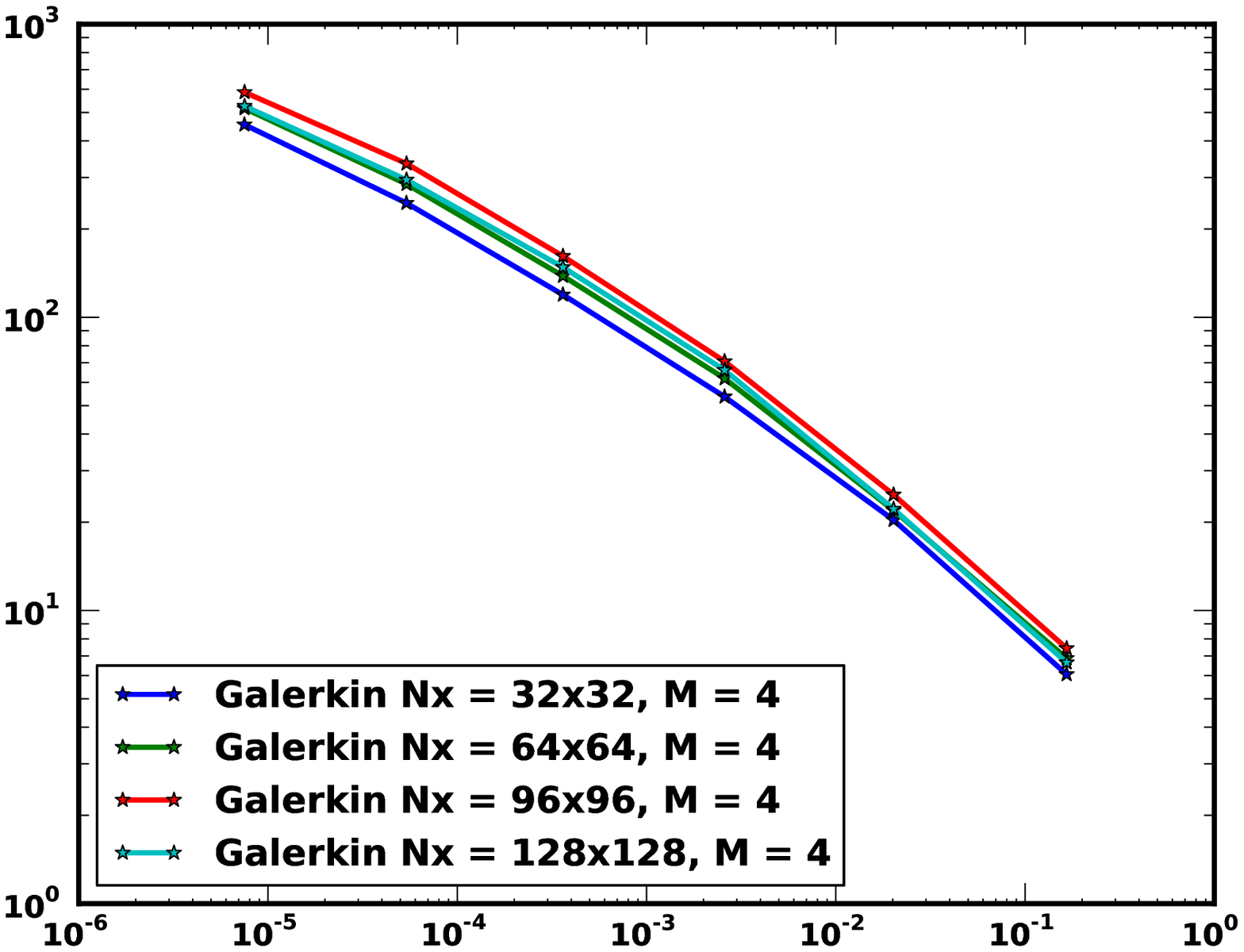}\label{RK:fig:Figure1b}
}
\caption{Scaled Galerkin solution time vs solution error for varying polynomial chaos order and spatial mesh for the diffusion coefficient modeled with uniform random variables}
\label{RK:fig:diff_unif_AGS_error} 
\end{center}
\end{figure}

For the numerical comparisons, we chose to discretize the domain $D=[-0.5,0.5]\times[-0.5,0.5]$ into a $64\times64$ grid resulting in a total number of nodes, $N_x = 4096$.  In figure~\eqref{RK:fig:diff_unif_AGS_error} the solution time for the stochastic Galerkin method, scaled by the deterministic solution time at the mean of the random field, is compared for different mesh sizes ($32\times32$, $64\times64$, $96\times96$ and $128\times128$) for the diffusion problem demonstrating that the solution time does not depend strongly on the mesh size (as is to be expected for the multi-grid preconditioner). The scaled solution time for these solvers and preconditioning techniques as a function of the standard deviation of the input random field, stochastic dimension, and polynomial order are then tabulated in Tables~\ref{RK:tab:time_dim_urf_ad}--\ref{RK:tab:time_sigma_lnrf_ad}.  In the tables, the number of Krylov iterations for the aforementioned preconditioners and iterations for Gauss-Seidel and Jacobi solvers are provided in parentheses. 
In the tables, MB, AGS, AJ,  GS and KP are the mean-based, approximate Gauss-Seidel, approximate Jacobi, Gauss-Seidel and Kronecker-product preconditioners respectively. GS in the solution methods refers to the Gauss-Seidel mean algorithm~(Algorithm-\ref{RK:alg:gs}).  ``Jacobi"  refers to the Jacobi mean algorithm~(Algorithm-\ref{RK:alg:jacobi}).  The solution tolerance for all of the stochastic Galerkin solvers is $1e^{-12}$.  For the Gauss-Seidel and Jacobi solvers, the inner solver tolerance is $3e^{-13}$.  All the computations are performed using a single core of an 8 core, Intel Xeon machine with 2.66 GHz and 16GB Memory.

\begin{figure}[!h]
\begin{center}
\includegraphics[scale=.65]{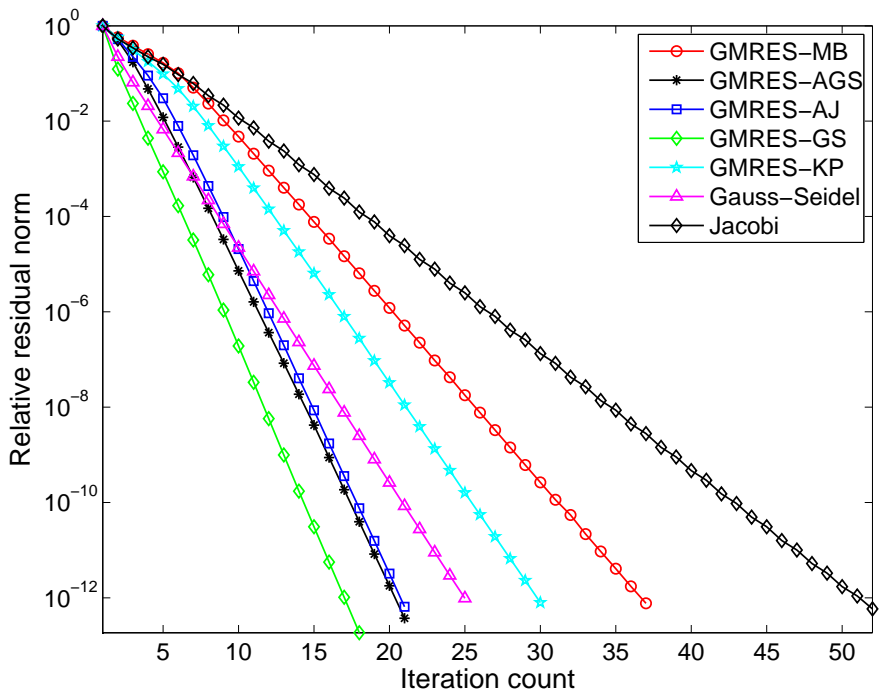} 
\caption{Relative residual norm vs iteration count for Galerkin system of equations, advection-diffusion equation with random field and uniform random variables} \label{RK:fig:ad_unif_iters}
\end{center}
\end{figure}

\begin{figure}[h]
\begin{center}
\includegraphics[scale=.65]{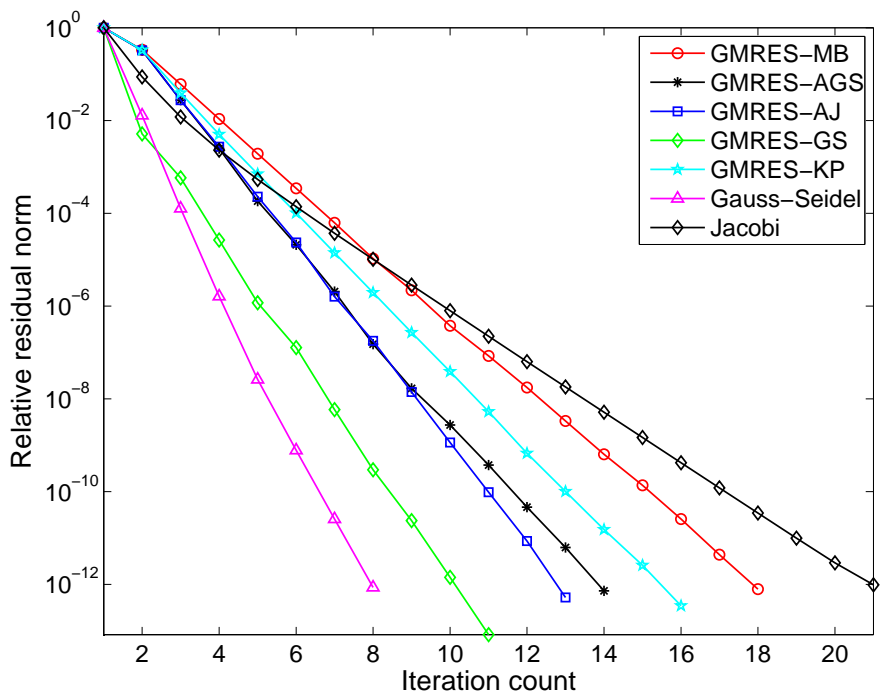} 
\caption{Relative residual norm vs iteration count for Galerkin system of equations, advection-diffusion equation with log-normal random field} \label{RK:fig:ad_logN_iters}
\end{center}
\end{figure}

Figures~\eqref{RK:fig:ad_unif_iters} and  \eqref{RK:fig:ad_logN_iters}  show the plots of relative residual error vs iteration count for the stochastic Galerkin system with stochastic dimension 4 and polynomial order 4 and standard deviation 0.1. It can be observed that the matrix-free Krylov solver with the Gauss-Seidel preconditioner takes the least number of iterations in case of uniform random field and Gauss-Seidel solver in case of log-normal random field, whereas the Jacobi solver takes highest number of iterations in both cases for a given tolerance. However in terms of solution time, the matrix-free Krylov solver with the approximate Gauss-Seidel preconditioner is the most efficient compared to all other stochastic Galerkin solvers. Comparison in terms of iteration count alone is misleading to evaluate the computational cost because, in each iteration  the cost of  preconditioner as observed in the case of Gauss-Seidel preconditioner could be very high resulting in higher computational cost even with small number of iterations. Hence we also compare the solution time for all preconditioners and solvers. 

In the tables, ``Div" means diverged.  In the case of diffusion coefficient modeled with uniform random random variables, with small variance ($\sigma=0.1$), it can be observed from Tables~\ref{RK:tab:time_dim_urf_ad} and~\ref{RK:tab:time_order_urf_ad} that more intrusive Krylov-based stochastic Galerkin solvers are more efficient than less intrusive Gauss-Seidel and Jacobi solvers.  Moreover the approximate Gauss-Seidel and Jacobi preconditioners are a significant improvement over the traditional mean-based approach.  However as the variance of the random field increases, we see from Table~\ref{RK:tab:time_sigma_urf_ad} that the Gauss-Seidel and Jacobi solvers suffer considerably, whereas the the Krylov-based approaches (excluding the Gauss-Seidel preconditioner) still perform quite well. This is not unexpected, as the operator becomes more indefinite as the variance increases.  Since the efficiency of the preconditioners and solvers based on Gauss-Seidel and Jacobi algorithms depends heavily on the efficiency of the deterministic solver, a good preconditioner for the mean stiffness matrix will significantly improve the above mentioned stochastic Galerkin preconditioners and solvers. 

In the case of the log-normal random field, we can see from Tables~\ref{RK:tab:time_dim_lnrf_ad} and~\ref{RK:tab:time_order_lnrf_ad} that Gauss-Seidel and Jacobi solvers have not performed well in terms of solution time. The Gauss-Seidel solver is comparable to the Krylov solver with the AGS preconditioner only in case of log-normal random field and at higher polynomial chaos order. For higher variance of the random field, we see from Table~\ref{RK:tab:time_sigma_lnrf_ad} that the Jacobi solver diverges. This problem might be addressed to some extent by using the true diagonal matrix $K^{k,k} = \sum_{i=0}^M c_{ikk} K_i$ from global stochastic stiffness matrix as the left-hand-side in the Jacobi solver and preconditioner instead of the mean matrix $K_0$.  From the Krylov iteration count provided in the parentheses of the tables, we can observe that the preconditioners and solvers based on Gauss-Seidel and Jacobi are robust with respect to stochastic dimension and polynomial order, however they are not robust with respect to the variance of the input random field.  The numerical results clearly show that the approximate Gauss-Seidel preconditioner performs betters than the mean-based preconditioner.  It was reported in \cite{RK:Rosseel2010} that the mean-based preconditioner was more efficient in terms of computational time than the Gauss-Seidel preconditioner. It could be due to the use of symmetric Gauss-Seidel preconditioner and duplicated matrix-vector products in the Gauss-Seidel preconditioner. In our work, we use non-symmetric Gauss-Seidel preconditioner which is half expensive as that of its symmetric version and we also minimize the duplicated matrix-vector products in Gauss-Seidel algorithm~(Algorithm-\ref{RK:alg:gs}).

\section{Conclusions}\label{RK:sec:conc}

In this work, various preconditioners for Krylov-based methods and solver methods based on Gauss-Seidel and Jacobi method are introduced. Results are compared with Krylov based methods (GMRES) with mean-based preconditioning.  The less intrusive Gauss-Seidel/Jacobi approaches did not perform well.  However these solvers can be used when legacy software has to be used to solve SPDEs with the stochastic Galerkin descritization. The use of approximate Gauss-Seidel or Jacobi preconditioners yields a significant improvement over traditional mean-based preconditioning. The Kronecker product preconditioner proposed in~\cite{RK:Ullmann2010} was in between mean-based preconditioner and approximate Gauss-Seidel preconditioner. Block and recycled Krylov methods present additional promising alternatives for improve the efficiency of the Jacobi and Gauss-Seidel solvers. 

\begin{center}
\bf\large Acknowledgements
\end{center}

This work was partially supported by the US Department of Energy through the NNSA Advanced Scientific Computing and Office of Science Advanced Scientific Computing Research programs.

\bibliographystyle{gCOM}

\clearpage

\begin{table}[h]
\caption{Scaled solution time ($\#$ of iterations) vs stochastic dimension for random field with uniform random variables, advection-diffusion equation, PC order = 4,  and $\sigma$=0.1}
\label{RK:tab:time_dim_urf_ad}
\begin{center}
\begin{tabular}{|c|c|c|c|c|c|c|c|c|}
 \hline
  Stoch. & \multicolumn{5}{|c|}{Preconditioners for GMRES} &
 \multicolumn{2}{|c|}{GS, Jacobi Solvers} \\ \cline{2-8} 
  dim & MB & AGS & AJ & GS & KP & GS &  Jacobi  \\ \hline
  1  & 12 (24)  & 7 (15)  & 14 (15)  & 65 (12)  & 9 (18)  & 64 (14)  & 114 (27)  \\ \hline
  2  & 47 (30)  & 27 (17)  & 48 (17)  & 240 (15)  & 37 (23)  & 256 (19)  & 489 (39)  \\ \hline
  3  & 133 (34)  & 74 (19)  & 129 (19)  & 631 (17)  & 106 (27)  & 727 (23) & 1353 (46)  \\ \hline
  4  & 291 (36)  & 165 (20)  & 285 (20)  & 1343 (18)  & 246 (29)  & 1579 (25)  & 3088 (52)  \\ \hline
  5  & 598 (38)  & 339 (21)  & 562 (21)  & 2571 (19)  & 532 (31)  & 3098 (27)  & 5941 (55)  \\ \hline
  6  & 1068 (40)  & 611 (22)  & 1012 (22)  & 4994 (21)  & 989 (34)  & 5808 (30)  & 11567 (62) \\ \hline
  7  & 1740 (41)  & 976 (22)  & 1690 (22)  & 7589 (21)  & 1619 (35)  & 9951 (32)  & 21199 (64)  \\
 \hline
\end{tabular}
\end{center}
\end{table}

\clearpage
\begin{table}[h]
\caption{Scaled solution time ($\#$ of iterations) vs order of polynomial chaos for random field with uniform random variables, advection-diffusion equation, Stoch. dim=4, $\sigma=0.1$}
\label{RK:tab:time_order_urf_ad}
\begin{center}
\begin{tabular}{|c|c|c|c|c|c|c|c|c|}
 \hline
  PC & \multicolumn{5}{|c|}{Preconditioners for GMRES} &
  \multicolumn{2}{|c|}{GS, Jacobi Solvers} \\ \cline{2-8}
  order & MB & AGS & AJ & GS & KP & GS & Jacobi  \\ \hline
  2  & 41 (26)  & 26 (16)  & 45 (16)  & 213 (13)  & 35 (22)  & 233 (17)  & 429 (33)  \\ \hline
  3  & 125 (32)  & 70 (18)  & 124 (18)  & 604 (16) & 102 (26)  & 663 (21)  & 1278 (43)  \\ \hline
  4  & 291 (36)  & 165 (20)  & 285 (20)  & 1343 (18)  & 246 (29)  & 1579 (25)  & 3088 (52)  \\ \hline
  5  & 626 (40)  & 336 (21)  & 591 (22)  & 2654 (20)  & 543 (32)  & 3302 (29)  & 6436 (60)  \\ \hline
  6  & 1163 (44)  & 637 (23)  & 1104 (24)  & 5082 (22)  & 1008 (35)  & 6311 (33)  & 12558 (68)  \\ \hline
  7  & 1981 (47)  & 1062 (24)  & 1825 (25)  & 8443 (23)  & 1713 (38)  & 11143 (36)  & 25117 (74)  \\
 \hline
\end{tabular}
\end{center}
\end{table}
\clearpage
\begin{table}[h]
\caption{Scaled solution time ($\#$ of iterations)  vs standard deviation ($\sigma$) for random field with uniform random variables, advection-diffusion equation, Stoch dim = 4, PC order = 4}
\label{RK:tab:time_sigma_urf_ad}
\begin{center}
\begin{tabular}{|c|c|c|c|c|c|c|c|c|}
 \hline
   & \multicolumn{5}{|c|}{Preconditioners for GMRES} &
 \multicolumn{2}{|c|}{GS, Jacobi Solvers} \\ \cline{2-8} 
  $\sigma$ & MB & AGS & AJ & GS & KP &GS & Jacobi  \\ \hline
  0.10  & 291 (36)  & 165 (20)  & 285 (20)  & 1343 (18)  & 246 (29)  & 1579 (25) & 3088 (52)  \\ \hline
  0.11  & 328 (41)  & 182 (22)  & 326 (23)  & 1556 (21)  & 284 (33)  & 1960 (31) & 3711 (62)  \\ \hline
  0.12  & 389 (49)  & 218 (26)  & 371 (26)  & 1921 (25)  & 325 (39)  & 2404 (38)  & 4624 (76)  \\ \hline
  0.13  & 469 (58)  & 257 (30)  & 461 (31)  & 2430 (30)  & 380 (46)  & 3039 (48)  & 5975 (97)  \\ \hline
  0.14  & 589 (73)  & 324 (38)  & 565 (39)  & 2968 (38)  & 478 (57)  & 4123 (65)  & 8153 (131)  \\ \hline
  0.15  & 810 (101)  & 437 (52)  & 742 (52)  & 3949 (52)  & 650 (78)  & 6110 (96)  & 12514 (196)  \\
 \hline
\end{tabular}
\end{center}
\end{table}

\clearpage

\begin{table}[h]
\caption{Sacled solution time ($\#$ of iterations) vs stochastic dimension for log-normal random field, advection-diffusion equation, PC order = 4 and $\sigma$=0.1}
\label{RK:tab:time_dim_lnrf_ad}
\begin{center}
\begin{tabular}{|c|c|c|c|c|c|c|c|c|}
 \hline
 Stoch. & \multicolumn{5}{|c|}{Preconditioners for GMRES} &
 \multicolumn{2}{|c|}{GS, Jacobi Solvers} \\ \cline{2-8}
  dim & MB & AGS & AJ & GS & KP &GS & Jacobi  \\ \hline
  1  & 9 (16)  & 8 (13)  & 12 (12)  & 48 (9)  & 12 (12)  & 36 (8)  & 81 (19) \\ \hline
  2  & 33 (17)  & 26 (13)  & 40 (12)  & 148 (9)  & 40 (12)  & 108 (8)  & 255 (20)  \\ \hline
  3  & 89 (17)  & 72 (13)  & 102 (12)  & 396 (10)  & 102 (12)  & 252 (8)  & 627 (21)  \\ \hline
  4  & 229 (17)  & 184 (13)  & 241 (12) & 858 (10)  & 215 (12)  & 507 (8)  & 1263 (21)  \\ \hline
  5  & 594 (17)  & 473 (13)  & 576 (12)  & 1765 (10)   & 572 (12)  & 915 (8)  & 2399 (22)  \\ \hline
  6  & 1426 (18)  & 1061 (13)  & 1221 (12)  & 3345 (10)  & 1370 (15)  & 1737 (8)  & 4044 (22)  \\
 \hline
\end{tabular}
\end{center}
\end{table}

\clearpage
\begin{table}[!h]
\caption{Scaled solution time ($\#$ of iterations)  vs order of polynomial chaos for log-normal random field,  advection-diffusion equation, Stoch. dim=4, $\sigma=0.1$}
\label{RK:tab:time_order_lnrf_ad}
\begin{center}
\begin{tabular}{|c|c|c|c|c|c|c|c|c|}
 \hline
 PC & \multicolumn{5}{|c|}{Preconditioners for GMRES} &
   \multicolumn{2}{|c|}{GS, Jacobi Solvers} \\ \cline{2-8}
  order & MB & AGS & AJ & GS & KP & GS & Jacobi  \\ \hline
  2  & 25 (15)  & 23 (13)  & 36 (12)  & 129 (8)  & 36 (12)  & 94 (7)  & 202 (16) \\ \hline
  3  & 77 (16)  & 65 (13)  & 97 (12)  & 353 (9) & 96 (12)  & 252 (8)  & 564 (19)  \\ \hline
  4  & 229 (17)  & 184 (13)  & 241 (12)  & 858 (10)  & 215 (12)  & 507 (8)  & 1263 (21)  \\ \hline
  5  & 733 (18)  & 585 (14)  & 645 (12)  & 1882 (10)  & 640 (12)  & 1030 (9)  & 2631 (24)  \\ \hline
  6  & 2157 (19)  & 1664 (14)  & 1784 (13)  & 4085 (10)  & 1938 (15)  & 1731 (9)  & 4821 (26)  \\
 \hline
\end{tabular}
\end{center}
\end{table}

\clearpage
\begin{table}[h]
\caption{Scaled solution time ($\#$ of iterations)  vs standard deviation ($\sigma$) for log-normal random field, advection-diffusion equation, Stoch dim = 4 and PC order = 4}
\label{RK:tab:time_sigma_lnrf_ad}
\begin{center}
\begin{tabular}{|c|c|c|c|c|c|c|c|c|}
 \hline
   & \multicolumn{5}{|c|}{Preconditioners for GMRES} &
 \multicolumn{2}{|c|}{GS, Jacobi Solvers} \\ \cline{2-8} 
  $\sigma$ & MB & AGS & AJ & GS & KP & GS & Jacobi  \\ \hline
  0.10  & 229 (17)  & 184 (13)  & 241(12)  & 858 (10)  & 215 (15)  & 507 (8)  & 1263 (21) \\ \hline
  0.15  & 275 (20)  & 214 (15)  & 267 (13)  & 1020 (12) & 266 (13)  & 634 (10)  & 2083 (34) \\\hline
  0.20  & 333 (24)  & 255 (18)  & 285 (14)  & 1177 (14)  & 285 (14)  & 825 (13)  & 3889 (62)  \\ \hline
  0.25  & 391  (28) & 285 (20)  & 305 (15)  & 1415 (17)  & 304 (15)  & 1018 (16)  & 17511 (254) \\ \hline
  0.30  & 464 (33)  & 328 (23)  & 361 (18)  & 1656 (20)  & 362 (18)  & 1210 (19)  & Div  \\ \hline
  0.35  & 529 (38)  & 387 (27)  & 660 (32)  & 1979 (24)  & 664 (32)  & 1464 (23)  & Div  \\
 \hline
\end{tabular}
\end{center}
\end{table}

\clearpage

\end{document}